\documentclass[10pt,a4paper,twoside,reqno]{elsart}

\usepackage[T1]{fontenc}
\usepackage[latin1]{inputenc}
\usepackage[english,francais]{babel}

\usepackage{lscape}     
\usepackage[french]{varioref} 

\usepackage{ifthen}      

\usepackage{enumerate}   

\usepackage{array}       
\usepackage{multirow}    
\usepackage{booktabs}

\usepackage{graphicx}    
\usepackage{rotating}   
\usepackage{color}
\usepackage{epsfig}
\usepackage{subfig}

\usepackage{float}

\usepackage{amsmath}     
\usepackage{amssymb}     
\usepackage{mathrsfs}    

\newtheorem{theoreme}{Th\'eor\`eme}[section]

\floatstyle{ruled}
\newfloat{algorithme}{hbt}{loalg}[section]
\floatname{algorithme}{\textbf{Algorithme}}

\newcommand{\guillemets}[1]{\og{}#1\fg{}} 

\newcommand{\abs}[1]{\ensuremath{\left|#1\right|}}
                                
\newcommand{\norme}[1]{\ensuremath{\left\|#1\right\|}} 

\newcommand{\dualite}[2]{\ensuremath{\left\langle#1,#2\right\rangle}}

\newcommand{\R}{\ensuremath{\mathbb{R}}}

\newcommand{\indi}[1]{\ensuremath{1\!\!1_{#1}}}
\renewcommand{\e}[1]{\ensuremath{\mathrm{e}^{#1}}}

\newcommand{\inc}{\ensuremath{\mathrm{inc}}}

\newcommand{\dif}{\ensuremath{\mathrm{dif}}}

\newcommand{\fonction}[5][]{\ifthenelse{\equal{#1}{}}%
  {\ensuremath{%
    \begin{array}[t]{@{}cll}
      #2&\longrightarrow&#3\\
      #4&\longmapsto&#5
    \end{array}}
  }
  {\ensuremath{%
    \begin{array}[t]{@{}lcll}
      #1\ :&#2&\longrightarrow&#3\\
      &#4&\longmapsto&#5
    \end{array}}
  }
}

\AtBeginDocument{}

\renewcommand{\phi}{\varphi}
\renewcommand{\epsilon}{\varepsilon}

\renewcommand{\O}[2][]{\ifthenelse{\equal{#1}{}}%
   {\ensuremath{\mathcal{O}\left(#2\right)}
   }
   {\ensuremath{\underset{#1}{\mathcal{O}}\left(#2\right)}
   }
}

\renewcommand{\o}[2][]{\ifthenelse{\equal{#1}{}}%
   {\ensuremath{\text{\scriptsize $\mathcal{O}$}\left(#2\right)}
   }
   {\ensuremath{\underset{#1}{\text{\scriptsize $\mathcal{O}$}}\left(#2\right)}
   }
}

\renewcommand{\i}{\ensuremath{\mathrm{i}}}

 {\everymath{\displaystyle\everymath{}}\array}%
 {\endarray}

\makeatletter
\renewcommand%
{\subsection}{\@startsection{subsection}{2}{0mm}
{-0.5\baselineskip}{-5mm}%
{\it\bfseries}}%
\makeatother

\journal{the Acad\'emie des sciences}
\begin{document}
\centerline{}
\begin{frontmatter}




%
\selectlanguage{francais}
\title{Imagerie laser}



\author[authorlabel1]{Jean-Baptiste Bellet},
\ead{bellet@cmap.polytechnique.fr}
\author[authorlabel2]{G\'erard Berginc}
\ead{gerard.berginc@fr.thalesgroup.com}

 \address[authorlabel1]{CMAP, \'Ecole Polytechnique, Route de Saclay, 91128 Palaiseau Cedex}
 \address[authorlabel2]{Thales Optronique, 2, Avenue Gay Lussac CS 90502, 78995 \'Elancourt Cedex}




\begin{abstract}
\selectlanguage{francais}
Nous \'elaborons une m\'ethode originale d'imagerie directe d'un objet. Il s'agit d'une m\'ethode de type d\'eriv\'ee topologique, destin\'ee \`a inverser des mesures contenant l'amplitude des ondes r\'etrodiffus\'ees \`a la suite d'illuminations laser.
\vskip 0.5\baselineskip

\selectlanguage{english}
\noindent{\bf Abstract}
\vskip 0.5\baselineskip
\noindent
{\bf Laser imaging. }
We derive an original direct imaging method of an object. It is based on topological derivatives, and aims at inverting the amplitude of  waves that are retropropagated after laser illuminations.
\end{abstract}


\end{frontmatter}

\selectlanguage{francais}

\section*{Introduction}

L'identification d'objets par un syst\`eme d'imagerie laser est un sujet aux domaines d'application vari\'es, comme la d\'etection d'engins explosifs improvis\'es (\emph{IED}) ou de tumeurs dans le derme. Dans ce contexte motivant, nous \'elaborons une m\'ethode d'imagerie directe (non it\'erative) par un tel syst\`eme. Ainsi, nous d\'eveloppons un algorithme d'inversion de mesures d'ondes r\'etrodiffus\'ees par l'objet \`a identifier, \`a la suite d'illuminations laser. Cette technique est monofr\'equentielle : une seule fr\'equence est exploit\'ee ; elle est monostatique : les sources d'illumination et les r\'ecepteurs sont confondus. Comme en pratique la vue est  partielle (une partie de l'objet n'est pas illumin\'ee), nous supposerons que les mesures sont prises sur une droite. Enfin, comme pour les applications, nous supposerons que les mesures contiennent l'amplitude (ou l'intensit\'e), mais ne contiennent pas la phase. Comme les m\'ethodes d'inversion directes usuelles (MUSIC \cite{AMM.INT}, backpropagation method, Kirchhoff migration) reposent sur la connaissance de la phase, le principal int\'er\^et et l'originalit\'e de notre m\'ethode concerne ce point.

Pour ajouter de l'information, de fa\c{c}on non restrictive par rapport aux applications, nous supposons que l'on peut prendre les mesures sur une autre droite, parall\`ele \`a la premi\`ere. Ainsi, pour chaque angle d'\'eclairage, nous disposons de deux mesures d'amplitude prises dans la m\^eme direction. En approximation haute fr\'equence, le passage d'une mesure d'amplitude \`a l'autre s'obtient par une \'equation de transport ; ceci cache une certaine notion de phase  dans le rapport entre les deux mesures. L'id\'ee est de trouver les petits disques, qui, \'eclair\'es, r\'etrodiffusent une onde dont la phase permet de \guillemets{transporter} les premi\`eres mesures vers les deuxi\`emes.  En plus d'utiliser directement les mesures, cette id\'ee permet d'exploiter l'information cach\'ee au sujet de la phase. Cette m\'ethode peut \^etre formul\'ee \`a l'aide d'un co\^ut que l'on cherche \`a rendre petit. La \guillemets{validit\'e} d'un point est alors test\'e par la valeur de la d\'eriv\'ee topologique du co\^ut : un point est bien plac\'e par rapport \`a l'objet recherch\'e si cette d\'eriv\'ee est tr\`es n\'egative. Il s'agit ainsi d'une approche par analyse asymptotique topologique, dans le m\^eme esprit que celle conduite dans \cite{AUR.IMA,AMM.IMA}.


\section{Position du probl\`eme}

\subsection{Probl\`eme direct.}
On consid\'ere l'espace libre $\R^2$ avec permittivit\'e $\epsilon_0>0$ et perm\'eabilit\'e $\mu_0>0$. On ajoute dans ce milieu une inclusion $D$ de permittivit\'e $\epsilon_*>0$ et perm\'eabilit\'e $\mu_*>0$. $D$ est suppos\'ee pr\`es de l'origine. On \'eclaire $D$ par un laser, \emph{i.e.}, une onde incidente plane d'angle unitaire $\eta$, de pulsation $\omega>0$ : $u_\inc(\eta,x)=\e{\i k_0\eta\cdot x}$. Le nombre d'onde associ\'e est $k_0:=\omega\sqrt{\epsilon_0\mu_0}$, et la longueur d'onde du laser est $\lambda:=\frac{2\pi}{k_0}$. Le nombre d'onde associ\'e \`a $D$ est $k_*:=\omega\sqrt{\epsilon_*\mu_*}$. L'onde qui se propage est alors $u[D](\eta,x)=u(x)$, solution de :
\begin{equation*}
    \nabla\cdot\left(\frac{1}{\mu_*}\indi{D}+\frac{1}{\mu_0}\indi{\R^2\setminus\overline{D}}\right)\nabla u+\omega^2(\epsilon_*\indi{D}+\epsilon_0\indi{\R^2\setminus\overline{D}})u=0,\\
\end{equation*}
avec une condition de radiation sortante sur l'onde diffract\'ee $u_\dif[D](\eta,x):=u[D](\eta,x)-u_\inc(\eta,x),x\in\R^2\setminus \overline{D}$. Num\'eriquement, ce probl\`eme est r\'esolu par la m\'ethode des \'equations int\'egrales, comme sugg\'er\'e dans \cite{BEL.MOD}.

\subsection{Probl\`eme inverse.}
On fixe deux droites de mesures parall\`eles en champ lointain : $\Gamma_0:x_2=\gamma_0$ et $\Gamma_1:x_2=\gamma_1$, o\`u les constantes $\gamma_0,\gamma_1>0$ sont grandes par rapport \`a la longueur d'onde, avec $\gamma_0-\gamma_1$ raisonnable (ni trop grand, ni trop petit). En \'eclairant l'objet $D$ avec l'angle $x/\abs{x}$, l'amplitude r\'etrodiffus\'ee en $x$ est $A[D](x)=\abs{u_\dif[D](x/\abs{x},x)}$, not\'ee $A_0[D]$ et $A_1[D]$ respectivement sur $\Gamma_0$ et $\Gamma_1$. On suppose qu'un objet $D_*$ est situ\'e \`a proximit\'e de l'origine, avec pour param\`etres \'electromagn\'etiques $\epsilon_*>0$ et $\mu_*>0$ connus. On suppose que l'on mesure $A_0[D_*]$ et $A_1[D_*]$. Le probl\`eme inverse que l'on souhaite r\'esoudre consiste \`a retrouver l'objet $D_*$ \`a partir des mesures, en inversant $A[D_*]$ par $A$, \emph{i.e.} en r\'esolvant $A[D]=A[D_*],$ sur $\Gamma_0\cup\Gamma_1$.


\section{M\'ethode de d\'eriv\'ee topologique}

\subsection{D\'emarche.}
Par approximation haute fr\'equence, on approche l'amplitude $A_1[D_*]$  en  transportant l'amplitude $A_0[D_*]$, avec une certaine phase $\phi[D_*]$ : si $A$ satisfait l'\'equation de transport : $\nabla A\cdot \nabla \phi[D_*]+\frac{1}{2}A\Delta\phi[D_*]=0$, avec  $A=A_0[D_*]$ sur $\Gamma_0$, alors la trace de $A$ sur $\Gamma_1$ approche $A_1[D_*]$. D'autre part, en transportant $A_0[D_*]$, avec la phase connue $\phi[D_\epsilon(z)]$  issue de l'onde diffract\'ee par un disque $D_\epsilon(z)$ de centre $z$ et de petit rayon $\epsilon>0$ : $\nabla A\cdot \nabla \phi[D_\epsilon(z)]+\frac{1}{2}A\Delta\phi[D_\epsilon(z)]=0$, avec $A=A_0[D_*]$ sur $\Gamma_0$, alors on obtient sur $\Gamma_1$ une amplitude $A_1[D_\epsilon(z)]$. On a alors l'intuition que $z$ est un point bien positionn\'e par rapport \`a l'objet $D_*$ lorsque l'amplitude transport\'ee $A_1[D_\epsilon(z)]$ est proche de l'amplitude mesur\'ee $A_1[D_*]$. L'id\'ee originale de reconstruction que nous proposons consiste \`a chercher les $z$ qui rendent petite la quantit\'e $\norme{A_1[D_\epsilon(z)]-A_1[D_*]}$. Nous obtenons une expression asymptotique de $A_1[D_\epsilon(z)]$ par la m\'ethode des caract\'eristiques. Puis un d\'eveloppement asymptotique de $\norme{A_1[D_\epsilon(z)]-A_1[D_*]}$ permet d'identifier une d\'eriv\'ee topologique \`a rendre tr\`es n\'egative pour estimer $D_*$.

\subsection{D\'eriv\'ee topologique.}
La r\'esolution de l'\'equation de transport par la m\'ethode des caract\'eristiques conduit au th\'eor\`eme suivant.
\begin{theoreme}
  Pour tout $p\in\Gamma_0$, en posant $\eta=\frac{z-p}{\abs{z-p}}$ l'angle d'\'eclairage, le transport de la mesure $A_0[D_*](p)$ vers $\Gamma_1$, le long de la direction de r\'etropropagation $(z,\eta)$, avec phase $\phi[D_\epsilon(z)](\eta,\cdot)$, conduit \`a l'amplitude :
\begin{equation*}
  A_1[D_\epsilon(z)](x_\eta)=A_0[D_*](p)\left(\frac{\gamma_0-z_2}{\gamma_1-z_2}\right)^{1/2}\left(1-\epsilon\frac{t_\eta}{4\abs{x_\eta-z}\abs{p_\eta-z}}\right)+\O{\epsilon^2}
\end{equation*}
o\`u $x_\eta=z+\frac{\gamma_1-z_2}{\eta_2}\eta=\Gamma_1\cap (z,\eta)$, et $t_\eta=-\frac{\gamma_1-\gamma_0}{\eta_2}$ est le temps de trajet.
\end{theoreme}
Param\'etrons par $x=(x_1,\gamma_1)$ : $ \eta(x)=\frac{z-x}{\abs{z-x}}$, $p(x)=z+\frac{\gamma_0-z_2}{\gamma_1-z_2}(x-z)$, et posons 
$  a_0(x,z):=A_0[D_*](p(x))\left(\frac{\gamma_0-z_2}{\gamma_1-z_2}\right)^{1/2}$, $ a_1(x,z):=-\frac{1}{4} a_0(x,z)\frac{\gamma_1-\gamma_0}{\gamma_0-z_2}\frac{1}{\abs{x-z}}$. Pour un domaine $D$,  notons $A_1[D]$ l'amplitude issue du transport de $A_0[D_*]$ sur $\Gamma_1$, dans la direction de r\'etropropagation, avec la phase $\phi[D]$ de l'onde diffract\'ee par $D$. Le th\'eor\`eme pr\'ec\'edent justifie le  r\'esultat asymptotique suivant :
\begin{theoreme}[D\'eriv\'ee topologique]
 La d\'eriv\'ee topologique de l'application $D\longmapsto \norme{A_1[D]-A_1[D_*]}$, \'evalu\'ee en $z$, dans la direction le disque unit\'e, est :
\begin{equation*}
  d(z)=\frac{1}{\norme{a_0(\cdot,z)-A_1[D_*]}}
\dualite{a_0(\cdot,z)-A_1[D_*]}{a_1(\cdot,z)}.
\end{equation*}
\end{theoreme}

\subsection{M\'ethode d'inversion directe.}
Si le petit disque $D_\epsilon(z)$ est bien plac\'e, \emph{i.e.} \guillemets{pr\`es} de $D_*$ alors la phase $\phi[D_\epsilon(z)]$ approche la phase $\phi[D_*]$, et donc $A_1[D_\epsilon(z)]$ approche $A_1[D_*]$. L'approche originale que nous proposons pour reconstruire la forme $D_*$ consiste  \`a minimiser $D\mapsto \norme{A_1[D]-A_1[D_*]}$, parmi les petits disques. Le d\'eveloppement asymptotique
du th\'eor\`eme : $$\norme{A_1[D_\epsilon(z)]-A_1[D_*]}=\norme{a_0(\cdot,z)-A_1[D_*]}+\epsilon d(z)+\O{\epsilon^2},$$ sugg\`ere que $\norme{A_1[D_\epsilon(z)]-A_1[D_*]}$ est petit lorsque la d\'eriv\'ee topologique $d(z)$ est n\'egative et tr\`es grande en valeur absolue. Pour estimer $D_*$, on trace ainsi les lignes de niveau de $z\longmapsto d(z)$, puis on en choisit une de niveau $d_0$ arbitraire. Celle-ci peut servir de point d\'epart pour un algorithme it\'eratif.



\section*{Remerciements}
Les auteurs remercient tr\`es chaleureusement Habib Ammari pour ses conseils et suggestions.

\bibliographystyle{elsart-num-sort}
\bibliography{biblio}
\addcontentsline{toc}{section}{R\'ef\'erences}

\end{document}